\documentstyle[12pt]{article}

\newcommand{\sect}[1]{\setcounter{equation}{0}\section{#1}}
\newcommand{\subsect}[1]{\subsection{#1}}

\renewcommand{\theequation}{\arabic{section}.\arabic{equation}}

\def\be{\begin{equation}}
\def\ee{\end{equation}}
\def\bea{\begin{eqnarray}}
\def\eea{\end{eqnarray}}

\def\1{\'{\i}}                           
\def\R{{\rm I\kern-.2em R}}

\def\a{\alpha}
\def\da{d_\alpha}

\def\ap{A_+}
\def\am{A_-}
\def\bp{B_+}
\def\bm{B_-}
\def\aa{N}
\def\bb{M}

\def\aap{a_+}
\def\aam{a_-}

\def\sch{{{\cal S}(1+1)}}

\begin{document}

\thispagestyle{empty}

 \
\hfill q-alg/9706014

\ 
\vspace{2cm}

\begin{center}
{\LARGE{\bf{ A Jordanian quantum }}}

{\LARGE{\bf{two-photon/Schr\"odinger algebra}}} 
\end{center}

\bigskip\bigskip

\begin{center} Angel Ballesteros$^\dagger$, Francisco J.
Herranz$^\dagger$ and Preeti Parashar$^\ddagger$
\end{center}

\begin{center} {\it {  $^\dagger$ Departamento de F\1sica, Universidad
de Burgos} \\   Pza. Misael Ba\~nuelos, 
E-09001-Burgos, Spain}
\end{center}

\begin{center} {\it {  $^\ddagger$ Department of Mathematics, Texas
A \& M University} \\ College Station, Texas 77843, 
USA}
\end{center}

\bigskip\bigskip\bigskip

\begin{abstract}
A non-standard quantum deformation of the two-photon algebra $h_6$ is 
constructed, and its quantum universal $R$-matrix is given.
Representations of this new quantum algebra are studied on the Fock
space and translated into Fock--Bargmann realizations that provide a
direct formalism for the definition of deformed states of light.
Finally, the isomorphism between
$h_6$ and the (1+1) Schr\"odinger algebra is used to introduce a new
(non-standard) Hopf algebra deformation of this   latter  symmetry
algebra.
\end{abstract}

\bigskip\bigskip\bigskip


\sect{Introduction}

The single-mode radiation field Hamiltonian that describes the
generation of a squeezed coherent state is given by
\cite{Gil}: \be
H=\hbar\,\omega (a_+\, a_- + \frac 1 2) + f_1(t)\,a_+ +
f_1(t)^\ast\,a_- + f_2(t)\,a_+^2 + f_2(t)^\ast\,a_-^2,
\label{ham}  
\ee
where $a_-$ and $a_+$ are the generators of the  boson algebra:
$[a_-,a_+]=1$. The dynamical symmetry algebra of this
Hamiltonian involving two-photon processes is the  ``two-photon algebra"
$h_6$, a six dimensional non-semisimple Lie algebra generated by
\be
\begin{array}{lll}
\aa=\aap\aam\qquad &\ap=\aap\qquad &\am=\aam\cr
\bb=1\qquad &\bp=\aap^2\qquad &\bm=\aam^2 \cr
\end{array}
\label{da}
\ee
where $\aa$ is the number operator and $\bb$ is a central
generator. The commutation rules among these generators are
\be
\begin{array}{lll}
 [\aa,\ap]=\ap \quad  &[\aa,\am]=-\am \quad &[\am,\ap]=\bb \cr
 [\aa,\bp]=2\bp \quad  &[\aa,\bm]=-2\bm \quad &[\bm,\bp]=4\aa+2\bb
\cr
 [\ap,\bm]=-2\am \quad  &[\ap,\bp]=0 \quad
&[\bb,\,\cdot\,]=0\cr
 [\am,\bp]=2\ap \quad  &[\am,\bm]=0 .\quad
& \cr
\end{array}
\label{ba}
\ee
 Three relevant subalgebras of $h_4$ arise:
the Heisenberg--Weyl algebra $h_3$ (generated by $A_+,A_-,M$), the
oscillator algebra $h_4$ (which is $h_3$ enlarged with the number
generator $N$) and the $su(1,1)$ algebra (defined by $B_+',B_-',N'$,
where $N'=(N+M/2)/2,\,{B_\pm}'={B_\pm}/2$). Thus  we have the
sequence $h_3\subset h_4\subset h_6$.  The essential role played by these
three symmetries in the algebraic description of coherent, squeezed and
intelligent states of light is well-known, and it has been recently
unified within a two-photon algebra approach in \cite{Brif}. Our
paper introduces a quantum deformation of the
two photon algebra and presents its potential algebraic abilities in
order to construct deformed (and non-classical) states of light.

In section 2, a non-standard quantum (Hopf algebra) deformation of
$h_6$ is obtained by making use of the Jordanian deformation of its
subalgebra $h_4$ \cite{osc}.  In this example, it will become clear
that although $B_+$ coincides with the square of $A_+$ in the
representation (\ref{da}), this relation is no longer true in general.
In spite of this fact, it is also shown that the quantum universal
$R$-matrix associated to the original $h_4$ deformation is just the
quantum
$R$-matrix for the two-photon quantization. 

In general, applications in Quantum Optics of all the
algebras   mentioned before, came initially from representations like
(\ref{da}) in terms of creation and   annihilation  operators (boson
realizations). Section 3 deals with the corresponding
deformation of the boson realization (\ref{da}). From an algebraic
point of view, the resultant representation can be translated into a
realization in terms of operators acting on the Hilbert space of
entire analytic functions where the construction of \cite{Brif} holds.
In this context, it is shown how the linear differential
operators of the undeformed case can be transformed after quantization
either into more complicated differential operators or into
differential-difference ones. The different role played
by $A_+$ and $A_-$ after quantization is emphasized, and the physical
 inequivalence  of algebraically equivalent quantizations is stressed.

 The  second main
purpose of the work is to show the usefulness (that, to our
knowledge, has not been used yet) of the isomorphism between the
two-photon algebra and the (1+1) Schr\"odinger algebra. In
particular, we shall present in the last section a new quantum
deformation of the latter algebra endowed with a Hopf algebra
structure, something that is left   out  in the previous
$q$-deformations so far obtained \cite{Vinet,Dobrev}. As a byproduct,
this isomorphism shows that the extended (1+1) Galilei algebra is also a
relevant subalgebra of $h_6$.


\sect{The deformation}

The non-standard (or Jordanian) quantum deformation of the
oscillator algebra $h_4$ was introduced in \cite{osc} (the
essentials of Lie bialgebra quantizations needed in order to derive
the following results can be found there; a comprehensive
presentation of the subject is given in \cite{PC}). The starting point
of the construction was the (coboundary) Lie bialgebra
$(h_4,\delta(r))$ defined by the classical $r$-matrix
\be
r=z\aa\wedge\ap \label{bb}
\ee
(which
is a solution of the classical Yang--Baxter equation) where $\aa$ and
$\ap$ belong to $h_4$. Therefore, the embedding $h_4\subset h_6$ allows
us to consider (\ref{bb}) as the generating object for a Jordanian
deformation of $h_6$ by taking $\aa$ and $\ap$ as generators for this
algebra. The cocommutators of the two-photon Lie bialgebra
$(h_6,\delta(r))$
are obtained from the relation
$\delta(X)=[1\otimes X+X\otimes 1,r]$ and read
\be
\begin{array}{l}
\delta(\ap)=0\qquad \delta(\bb)=0\cr
\delta(\aa)=z\aa\wedge\ap\qquad 
\delta(\bp)=-2z\bp\wedge\ap\cr
\delta(\am)=z(\am\wedge\ap + \aa\wedge\bb)\cr
\delta(\bm)=2z(\bm\wedge\ap + \aa\wedge\am).\cr
\end{array}
\label{bc}
\ee

Note that $\delta(\bm)$ contains the term $\aa\wedge\am$, which
involves two ``non-primitive" generators. As a consequence (see
\cite{osc}) it is not possible to ``exponentiate" directly the first
order in the deformation given by the Lie bialgebra (\ref{bc}) in
order to obtain the full coproduct of the quantum deformation we are
looking for. However, by assuming that the right quantum structure is
not far away from the exponentiation of (\ref{bc}), the following
coassociative coproduct for the quantum two-photon algebra $U_z(h_6)$
can be obtained:
\be
\begin{array}{l}
\Delta(\ap)=1\otimes \ap + \ap \otimes 1\qquad
\Delta(\bb)=1\otimes \bb + \bb \otimes 1\cr
\Delta(\aa)=1\otimes \aa + \aa \otimes e^{z\ap}
\qquad
\Delta(\bp)=1\otimes \bp + \bp \otimes e^{-2z\ap}\cr
\Delta(\am)=1\otimes \am + \am \otimes e^{z\ap}
+ z\aa \otimes e^{z\ap} \bb\cr
\Delta(\bm)=1\otimes \bm + \bm \otimes e^{2z\ap}
+ z\aa \otimes e^{z\ap}(\am-z \bb\aa)
-z\am\otimes e^{z\ap} \aa .\cr
\end{array}
\label{bd}
\ee
Deformed commutation rules consistent with (\ref{bd}) are given in
the form:
\bea
&& [\aa,\ap]=\frac {e^{z\ap}-1}z \qquad   [\aa,\am]=-\am \qquad
 [\am,\ap]=\bb  e^{z\ap} \cr
&& [\aa,\bp]=2\bp \qquad   [\aa,\bm]=-2\bm - z \am \aa
\qquad [\bb,\,\cdot\,]=0\cr
&&[\bm,\bp]= 2(1+  e^{-z\ap})\aa + 2\bb - 2 z \am\bp  \label{bg}\\
&& [\ap,\bm]=-(1+  e^{z\ap}) \am + z e^{z\ap}\bb\aa\qquad 
[\ap,\bp]=0\cr
&&  [\am,\bp]=2\frac {1-e^{-z\ap}}z 
  \qquad   [\am,\bm]=-z\am^2 . 
\nonumber 
\eea
Finally,  counit $\epsilon$ and  antipode $\gamma$ for $U_z(h_6)$ read
\be
\epsilon(X) =0 \qquad 
 \mbox{for\ \ $X\in \{\aa,\ap,\am,\bp,\bm,\bb\}$} 
\label{be}
 \ee
\be
\begin{array}{l}
\gamma(\ap)=-\ap\qquad \gamma(\bb)=-\bb\cr
\gamma(\aa)=-\aa e^{-z\ap}\qquad
\gamma(\bp)=-\bp e^{2z\ap}\cr
\gamma(\am)=-\am e^{-z\ap}+z\aa\bb e^{-z\ap}\cr
\gamma(\bm)=-\bm e^{-2z\ap}-z\am e^{-2z\ap} .\cr
\end{array}
\label{bf}
\ee

Note that the non-standard quantum oscillator algebra $U_z(h_4)$
is a Hopf subalgebra of $U_z(h_6)$. The quantum universal $R$-matrix for
$U_z(h_4)$
\be
{\cal R}=\exp\{-z\ap\otimes\aa\}\exp\{z\aa\otimes\ap\}  
\label{ca}
\ee
was found in \cite{osc}. In fact, this element ${\cal R}$ can be proven
to be the quantum $R$-matrix for $U_z(h_6)$ when $\ap$ and $\aa$ are
considered as generators of the latter algebra (see Appendix A). 
Recall that its classical counterpart (\ref{bb}) has been used to
define the whole Lie bialgebra (\ref{bc}), and that this kind of
 embedding  has been
shown to preserve the quantum $R$-matrix for the initial subalgebra
in other examples of non-standard deformations (see
\cite{Tmatrix,Rmatrix}).


\sect{Representation theory}

The physical relevance of the two-photon algebra comes from the
Hamiltonians of the type (\ref{ham}), in which a certain
boson representation for $h_6$ gives a dynamical symmetry 
of the problem. Let us generalize such kind of representations to the
quantum case.

\subsect{Quantum boson realizations}

 A  realization of   $h_6$ (\ref{ba})  in terms of the boson generators
$\aam$, $\aap$ is given by (\ref{da}). When the previous boson operators
act on the number states Hilbert space spanned by
$\{|m\rangle\}_{m=0}^{\infty}$ the action of the two-photon operators on
these states is
\be
\begin{array}{l}
\aa |m\rangle =m\, |m\rangle \qquad
\bb |m\rangle =  |m\rangle\cr
\ap |m\rangle = \sqrt{m+1}\,|m+1\rangle  \qquad 
\am |m\rangle = \sqrt{m}\,|m-1\rangle  \cr
\bp |m\rangle = \sqrt{(m+1)(m+2)}\,|m+2\rangle  \qquad 
\bm |m\rangle = \sqrt{m(m-1)}\,|m-2\rangle  .\cr
\end{array}
\label{dc}
\ee

A part of this representation can be deformed by considering
the
results given in \cite{bosones} concerning the subalgebra defined
by the quantum oscillator  $U_z(h_4)$. In particular, by setting 
$\beta=0$, $\delta=1$, and by replacing $2z\to
z$ in those results, we get the realization: 
\be
 \aa= \frac{e^{z\aap}-1}{z}\,\aam \qquad
  \ap=\aap  \qquad 
  \am=  e^{z\aap}\aam  \qquad  \bb= 1 .
\label{dd}
\ee 
The remaining generators are easily found to be realized as
\be
\bp=\left(\frac{1- e^{-z\aap} }{z}\right)^2  \qquad 
\bm=e^{z\aap}\aam^2 .
\label{ddd}
\ee

Therefore, the classical identification $B_+=A_+^2$ is no longer valid
in the quantum case. Moreover,
the coproduct (\ref{bd}) for  $\Delta(B_+)$ is strongly different from
$\Delta(A_+)^2$ (in fact, this property appears at a non-deformed level
aswell \cite{Enrico}). 

Now, by taking
into account that 
 \be 
e^{ z \aap}|m\rangle =|m\rangle+   
\sum_{k=1}^\infty \frac{z^k }{k!} \, \sqrt{\frac {(m+k)!}{m!} }
\, |m+k\rangle  
\label{de}
\ee
 the action of the generators of $U_z(h_6)$ on the number states
$\{|m\rangle\}_{m=0}^{\infty}$  reads
\bea
&&\ap|m\rangle =\sqrt{m+1}\, |m+1\rangle  \qquad
\bb|m\rangle =  |m\rangle   \cr
&&\am|m\rangle = \sqrt{m}\,
|m-1\rangle + m \sum_{k=0}^\infty
 \frac{z^{k+1} }{(k+1)!} \, \sqrt{\frac {(m+k)!}{m!} }
 \, |m+k\rangle  \cr
&&\aa|m\rangle =m\, |m\rangle + m \sum_{k=1}^\infty  \frac{z^k
}{(k+1)!} \, \sqrt{\frac {(m+k)!}{m!} }
 \, |m+k\rangle  \cr
&&\bp|m\rangle =\sqrt{(m+1)(m+2)}\,|m+2\rangle \cr
&&\qquad \qquad + \sum_{k=1}^\infty
(-2+2^{k+2}) \frac{(-z)^k }{(k+2)!} \, \sqrt{\frac {(m+k+2)!}{m!} }
 \, |m+k+2\rangle  \cr
&&\bm|m\rangle =\sqrt{m (m-1)}\, |m-2\rangle +
z \sqrt{m} (m-1)\, |m-1\rangle \cr
&&\qquad \qquad + m (m-1) \sum_{k=0}^\infty  \frac{z^{k+2}
}{(k+2)!} \, \sqrt{\frac {(m+k)!}{m!} }
 \, |m+k\rangle  .\label{df}
\eea

The numbers $\langle m|X|m'\rangle$ where 
$\langle m|m'\rangle=\delta_{m,m'}$ give rise to the matrix  elements
of this infinite-dimensional representation; explicitly,
$$
 \ap=\left(\begin{array}{llllll}
0 & . & . & . & . & \cr 
1  & 0 & . & . & . &\cr . 
& {\sqrt{2}} & 0 & . & . & \cr
. & . & {\sqrt{3}} & 0 & . &\cr 
. & . & . &{\sqrt{4}} & 0 &\cr
  &    &   &  &   &\ddots\cr
\end{array}\right)\ 
 \bp=\left(\begin{array}{llllll}
0 & . & . & . & . & \cr
0  & 0 & . & . & . &\cr
{\sqrt{2}} & 0 & 0 & . & . & \cr
-{\sqrt{6}} z & {\sqrt{6}} & 0 & 0 & . & \cr
\frac 7 {\sqrt{6}} {z^2} & -2 {\sqrt{6}} z &
2 {\sqrt{3}} & 0 & 0 & \cr
  & & & & &\ddots\cr
\end{array}\right)
$$
$$
 \am=\left (\begin{array}{llllll}
.  & 1 & . & . & .& \cr 
. &   z  & \sqrt{2} & . & .& \cr 
. &
  \frac 1{\sqrt{2}}   z^2 &  2  
  z   & \sqrt{3} & . &\cr 
. &  \frac 1 {\sqrt{6}}     {z^3} & 
  {\sqrt{3}}        z^2 &  3   z 
    & \sqrt{4} &\cr 
. & 
   \frac 1  {\sqrt{24}}    
   {z^4} &    {2\over {\sqrt{3}}}   
    {z^3} & 
      3  z^2 &  4   z & \cr
  &    &   &  &   &\!\!\ddots\cr
\end{array} \right) \ 
\bm=\left(\begin{array}{llllll}
. & . & {\sqrt{2}} & . & . & \cr
. & . & {\sqrt{2}} z &
{\sqrt{6}} & . & \cr . & . & z^2 &
2 {\sqrt{3}} z & 2 {\sqrt{3}} & \cr
. & . & \frac 1{\sqrt{3}} z^3 & 3 z^2 & 6 z & \cr
. & . & \frac 1 {2 \sqrt{3}} z^4 & 2 z^3 & 6 z^2 & \cr
& & & & &\!\!\ddots\cr
\end{array} \right)
$$
\be
 \aa=\left (\begin{array}{llllll}
.  & . & . & . & .& \cr 
. &    1   & . & . & .& \cr 
. &
  \frac 1{\sqrt{2}}     z & 2  
    & . & . &\cr 
. & 
       {1\over {\sqrt{6}}}   {z^2} & 
  {\sqrt{3}}      z & 3  
    & . &\cr 
. & 
   \frac 1  {\sqrt{{24}}}     
   {z^3} & \frac 2 {\sqrt{3}} 
      {z^2} & 
       3  z & 4  & \cr
  &    &   &  &   & \ddots\cr
\end{array} \right)  \quad
 \bb=\left(\begin{array}{llllll}
1 & . & . & . & . & \cr 
.  & 1 & . & . & . &\cr . 
& . & 1 & . & . & \cr
. & . & . & 1 & . &\cr 
. & . & . &. & 1 &\cr
  &    &   &  &   &\ddots\cr
\end{array}\right).
\ee
Note the characteristic appearance of monomials in $z$ as entries of
the deformed matrices. In this respect, see
\cite{bosones,dob,ab,abc,aizawa,vander} for representations of
non-standard quantum algebras where this feature appears
recurrently.

\subsect{Fock--Bargmann realizations}

The deformed boson realization (\ref{dd}--\ref{ddd}) can be immediately
translated into differential operators acting on the the space of
entire analytic functions $f(\a)$ (the Fock--Bargmann
(FB) representation
\cite{FB}):
\bea
&& \aa= \frac{e^{z\a}-1}{z}\,\frac{d}{d\a}  \qquad
  \ap=\a   \qquad 
  \am=  e^{z\a}\,\frac{d}{d\a}   \qquad  \bb= 1  \cr
&& \bp=\left(\frac{1- e^{-z\a} }{z}\right)^2   \qquad 
\bm=e^{z\a}\,\frac{d^2}{d\a^2} .
\label{fbr}
\eea
These expressions would be the starting point for a general study of
deformed states of light generated by the quantum
deformation introduced here (see \cite{Brif} for the classical
construction). However, there also exists another set of operators
defined on the FB space and algebraically linked to this
deformation. It is easy to check that the following automorphism of
$h_6$
\be
\begin{array}{lll}
N\rightarrow - N&\qquad A_+\rightarrow - A_-&\qquad A_-\rightarrow
- A_+\cr
M\rightarrow - M&\qquad B_+\rightarrow - B_-&\qquad B_-\rightarrow
- B_+   \cr
\end{array}
\label{aut}
\ee
and the transformation of the deformation parameter
\be
z\to -z
\label{autb}
\ee
transforms the Lie bialgebra (\ref{bc}) into 
\be
\begin{array}{l}
\delta(\am)=0\qquad \delta(\bb)=0\cr
\delta(\aa)=z\aa\wedge\am\qquad 
\delta(\bm)=-2z\bm\wedge\am\cr
\delta(\ap)=z(\ap\wedge\am + \aa\wedge\bb)\cr
\delta(\bp)=2z(\bp\wedge\am + \aa\wedge\ap).\cr
\end{array}
\label{bcp}
\ee
 Therefore, (\ref{bcp}) and (\ref{bc}) are equivalent
as Lie bialgebras, and the quantization of the former leads to
expressions that can be deduced from the ones given in section 2 by
making use of the automorphism (\ref{aut}) and  (\ref{autb}). The
essential feature of this deformation is that now $A_-$ and $M$ are the
primitive generators. We shall write only the set of resultant deformed
commutation rules: 
\bea
&& [\aa,\am]=-\frac {e^{z \am}-1}{z} \qquad
[\aa,\ap]=\ap \qquad
[\am,\ap]=\bb  e^{z\am} \cr
&& [\aa,\bm]=-2\bm \qquad
[\aa,\bp]=2\bp + z \ap \aa
\qquad [\bb,\,\cdot\,]=0\cr
&&[\bp,\bm]=- 2(1+  e^{-z\am})\aa - 2\bb + 2 z \ap\bm  \label{bgam}\\
&& [\am,\bp]=(1+  e^{z\am}) \ap - z e^{z\am}\bb\aa\qquad 
[\am,\bm]=0\cr
&&  [\ap,\bm]=-2\frac {1-e^{-z\am}}{z} 
  \qquad   [\ap,\bp]=z\ap^2 . 
\nonumber 
\eea

The corresponding
FB realization for (\ref{bgam})  
reads (we write
$\frac{d}{d\a}=\da$): 
\bea
&& \aa= \a\,\left(\frac{e^{z\da}-1}{z}\right) 
\qquad
  \ap=\a\,e^{z\da}   \qquad 
  \am=  \da  \qquad  \bb= 1  \cr
&& \bp=(\a^2  + z\,\a)\,e^{z\da} 
\qquad 
\bm=\left(\frac{1-e^{-z\da}}{z}\right)^2 .
\label{fbram}
\eea
The deformed boson realization would be obtained by substituting
$\da=a_-$ and $\a=a_+$. Note that (\ref{fbram}) is a
differential-difference realization; in particular, the action
\be
\left(\frac{e^{z\da}-1}{z}\right)\,f(\a)=\frac{ f(\a+z)-f(\a)}{z},
\label{disc}
\ee
corresponds to a discrete derivative. Therefore, (\ref{fbram}) can
be thought   of  as a certain discretization of the usual FB
representation. 

From a physical point of view, relevant states of the radiation
field
appear as eigenfunctions of the generators of the two-photon algebra in
the FB representation \cite{Brif}. Therefore, the two (algebraically
equivalent) quantum deformations of
$h_6$  presented here  give rise to different deformed
eigenproblems on the space of analytic functions. This procedure seems
to be the most natural way for generating Jordanian (and, in general,
quantum algebra) analogues of coherent, squeezed and intelligent
states. A
detailed analysis of this construction and of its physical contents will
be presented elsewhere, but we can advance here that the former
(\ref{fbr}) ($\ap$ primitive) would originate a class of smooth deformed
states, and the latter (\ref{fbram}) ($\am$ primitive) will be linked to
a set of states including some intrinsic discretization. In both cases,
the precise shape of the deformation is a consequence of the
compatibility with the deformed composition rule of representations given
by the corresponding quantum coproduct.


\sect{A quantum (1+1) Schr\"odinger algebra}

The (1+1) dimensional Schr\"odinger algebra $\sch$
\cite{Schroda,Schrodb} is a six-dimensional Lie algebra generated by $H$
(time translation), $P$ (space translation),
$K$ (Galilean boost),
$D$ (dilation),
$C$ (conformal transformation) and $M$ (mass); it is endowed with the
following commutation rules:
\be
\begin{array}{lll}
 [D,P]=-P \quad  &[D,K]=K \quad &[K,P]=M \cr
 [D,H]=-2H \quad  &[D,C]=2C  \quad &[H,C]=D \cr
 [K,H]=P \quad  &[K,C]=0 \quad
&[M,\,\cdot\,]=0\cr
 [P,C]=-K \quad  &[P,H]=0 .\quad
& \cr
\end{array}
\label{fa}
\ee
The subalgebra generated by $\{ H,P,K,M\}$ defines a (1+1) extended
Galilei algebra and the one spanned by $\{ D,C,H\}$ closes as an
$sl(2,\R)$ structure. 

The physical interest  in $\sch$ comes from the fact that (\ref{fa})
is the Lie algebra of the symmetry group of the (1+1) dimensional
Schr\"odinger equation. Recently, the symmetry analysis of some
differential-difference generalizations of this outstanding equation has
originated various $q$-analogues of $\sch$, none of them endowed
with a known Hopf algebra structure \cite{Vinet,Dobrev}. However,   an
isomorphism  between $\sch$ and $h_6$
\cite{negro}, can be explicitly given as follows:
\be
D=-\aa-\frac 12\bb \qquad P=\ap\qquad K=\am\qquad
H=\frac 12\bp\qquad C=\frac 12\bm 
\label{fb}
\ee
keeping $M$ as the central generator for both algebras.
Therefore, starting from the results presented in previous sections
and by using (\ref{fb}), it is immediate to obtain a
complete (non-standard) quantum deformation of $\sch$. Its Lie
bialgebra is provided by the classical $r$-matrix
\be
r=z P\wedge D +\frac z2 P\wedge M
\label{fc}
\ee
with cocommutators given by
\bea
&&\delta(P)=0\qquad \delta(M)=0\qquad
\delta(H)=-2z H\wedge P\cr
&& \delta(K)=z(K\wedge P - D\wedge M)\qquad
\delta(D)=z(D\wedge P +\frac 12 M\wedge P)\cr
&&\delta(C)=z(2C\wedge P + K\wedge D + \frac 12 K\wedge M) .
\label{fd}
\eea
The coproduct and the commutation rules of the Hopf algebra
$U_z(\sch)$ are
\bea
&&\Delta(P)=1\otimes P + P \otimes 1
\qquad \Delta(M)=1\otimes M + M \otimes 1\cr
&&\Delta(H)=1\otimes H + H\otimes e^{-2zP}\cr
&&\Delta(K)=1\otimes K + K\otimes e^{zP}-z(D+\frac 12 M)\otimes
e^{zP}M\cr
&&\Delta(D)=1\otimes D + D\otimes e^{zP} + \frac 12 M\otimes
(e^{zP}-1)\cr
&&\Delta(C)=1\otimes C + C\otimes  e^{2zP}-\frac z2 (D + \frac 12
M)\otimes  e^{zP}(K+z(D + \frac 12 M)M)\cr
&&\qquad\quad +\frac z2 K\otimes  e^{zP}(D + \frac 12 M),
\label{fe}
\eea
\bea
&& [D,P]=\frac {1- e^{z P}}z \qquad   [D,K]=K\qquad
 [K,P]=M  e^{z P } \cr
&& [D,H]= - 2 H \qquad   [D,C]=2 C -\frac z2 K(D+\frac 12 M)  \cr
&&[H,C]= \frac 12(1+  e^{-z P})(D+\frac 12 M) -\frac 12 M +  z K
H\cr  
&& [K,H]=\frac {1-e^{-zP}}z \qquad  [K,C]=-\frac z2 K^2
\qquad [M,\,\cdot\,]=0\cr
&&  [P,C]=-\frac 12 (1+ e^{zP}) K - \frac z2 e^{zP} M(D+\frac 12 M)
  \qquad   [P,H]=0 . 
\label{ff}
\eea
 
Note that neither $\{ H,P,K,M\}$ (the extended Galilei algebra)
nor
$\{ D,C,H\}$ define a Hopf subalgebra of this deformation. However,
the generators
$\{D,P,K,M\}$ do close as a Hopf subalgebra (they were the
oscillator
generators in $h_6$). Finally, the (factorized) universal $R$-matrix
for
$U_z(\sch)$ is written as
\be
\begin{array}{l}
{\cal R}=\exp\{zP\otimes D\}\exp\{\frac z2 P\otimes M \}
\exp\{-\frac z2 M\otimes P \}\exp\{-z D\otimes P\} .\cr
\end{array}
\ee

This connection between the  two-photon  and the
(1+1) Schr\"odinger algebras deserves  further  analysis. From it,
the
extended (1+1) Galilei algebra appears as another relevant subalgebra
of $h_6$ and, in turn, the oscillator algebra $h_4$ can be embedded
within $\sch$. On the other hand, the isomorphism (\ref{fb}) 
exhibits again, in the geometrical language of the Schr\"odinger
symmetry, the physical   inequivalence  between
the (algebraically equivalent) quantizations (\ref{bg}) and (\ref{bgam})
shown in section 3 by constructing FB realizations: the first
deformation makes $\ap=P$ the primitive generator, a role played in the
second one by the boost $\am=K$.



\bigskip
\bigskip

\noindent
{\Large{{\bf Acknowledgments}}}

\bigskip

A.B. and F.J.H. have been
partially supported by DGICYT (Project  PB94--1115) from the
Ministerio de Educaci\'on y Cultura de Espa\~na and by Junta de
Castilla y Le\'on (Projects CO1/396 and CO2/297). P.P. wishes to thank 
Prof. Edward Letzter for hospitality at Texas A \& M University.

\bigskip
\bigskip

\noindent
{\Large{{\bf Appendix A}}}

\appendix

\setcounter{equation}{0}

\renewcommand{\theequation}{A.\arabic{equation}}

\bigskip

\noindent
As it was shown in  \cite{osc}  the element (\ref{ca})
  fulfills both the quantum Yang--Baxter equation and the
property
\be
 {\cal R}\Delta(X){\cal R}^{-1}=\sigma\circ \Delta(X) \qquad
\mbox{for}\ \ X\in\{\aa,\ap,\am,\bb\} 
\label{cb}
\ee
where $\sigma$ is the flip operator: $\sigma (a\otimes b)=b\otimes a$.
Since $U_z(h_4)\subset
U_z(h_6)$ we only need to prove  that  $\cal R$ also satisfies
(\ref{cb}) for the two remaining generators $\bp$ and $\bm$.
We shall make use of the formula
 \be
e^{f}\,\Delta(X)\,e^{-f}=\Delta(X) +\sum_{n=1}^\infty \frac
1{n!}\,[f,[\dots [f,\Delta(X)]]^{n)}\dots] 
\label{cj}
\ee
where the superindex $n)$ means $n$ commutators.
To begin with we take   $X\equiv \bp$ and  $f\equiv
z\aa\otimes\ap$. Straightforward computations lead to
\be
[z\aa\otimes\ap,[\dots [z\aa\otimes\ap,\Delta(\bp)]]^{n)}\dots]
=\bp\otimes  e^{-2z\ap}(2z\ap)^n\quad n\ge 1
\ee
so that 
\bea
&&e^{z\aa\otimes\ap}\Delta(\bp)e^{-z\aa\otimes\ap}
=\Delta(\bp)
+\sum_{n=1}^\infty \bp\otimes  e^{-2z\ap}\frac {(2z\ap)^n}{n!}\cr
 &&\qquad = 1\otimes \bp + \bp \otimes e^{-2z\ap}+
\bp \otimes e^{-2z\ap}(e^{2z\ap}-1)\cr
&&\qquad =1\otimes \bp + \bp \otimes
1\equiv \Delta_0(\bp) .\label{ck}
\eea
On the other hand, we find that
\be
[-z\ap\otimes\aa,[\dots [-z\ap\otimes\aa,\Delta_0(\bp)]]^{n)}\dots]
=(-2z\ap)^n  \otimes \bp \quad n\ge 1
\ee
and the proof for $\bp$ follows
\bea
&&e^{-z\ap\otimes\aa}\Delta_0(\bp)e^{z\ap\otimes\aa}
=\Delta_0(\bp)
+\sum_{n=1}^\infty \frac {(-2z\ap)^n}{n!} \otimes \bp\cr
 &&\qquad = 1\otimes \bp + \bp \otimes1 +
(e^{-2z\ap}-1)  \otimes \bp =\sigma\circ \Delta(\bp)  .\label{cl}
\eea

We consider now  $X\equiv \bm$  and  $f\equiv
z\aa\otimes\ap$; thus we have  
\be
\begin{array}{l}
[z\aa\otimes\ap, \Delta(\bm)]
=z\aa\otimes \{ e^{z\ap}(zMN-\am)-\am\} - z^2\aa^2\otimes
e^{z\ap} M\cr
\qquad\qquad   -2z\bm\otimes  e^{2z\ap} \ap 
 -z^2\am\aa\otimes e^{2z\ap} \ap\cr
\qquad\qquad   +z\am \aa\otimes ( e^{2z\ap} -  e^{z\ap})
+z^2\am\otimes e^{z\ap} \ap \aa ,\cr
\end{array} 
\label{cm}
\ee 
\be
\begin{array}{l}
  [z\aa\otimes\ap,[z\aa\otimes\ap, \Delta(\bm)]]
=2 z^2\aa^2\otimes
e^{z\ap} M+ 4z^2\bm\otimes  e^{2z\ap} \ap^2 \cr
   + 3 z^3\am\aa\otimes e^{2z\ap} \ap^2 
  - 2 z^2\am \aa\otimes ( e^{2z\ap} -  e^{z\ap})\ap 
   -z^3\am\otimes e^{z\ap} \ap^2 \aa .\cr
\end{array} 
\label{cn}
\ee 
In general, a recurrence method gives (for $n\ge 3$):
\be
\begin{array}{l}
 [z\aa\otimes\ap,[\dots [z\aa\otimes\ap,\Delta(\bm)]]^{n)}\dots]
= \bm\otimes  e^{2z\ap} (-2z\ap)^n \cr
  + z(2^n-1)  \am\aa\otimes e^{2z\ap}
(-z\ap)^n +
 z n\am \aa\otimes ( e^{2z\ap} -  e^{z\ap})(-z\ap)^{n-1}\cr
  -z \am\otimes e^{z\ap} (-z\ap)^n \aa .
\end{array} 
\label{co}
\ee 
 From (\ref{cm}--\ref{co}), the following result is obtained:
\bea
&&g=e^{z\aa\otimes\ap}\,\Delta(\bm)\,e^{-z\aa\otimes\ap}
=\Delta(\bm) + z\aa\otimes \{ e^{z\ap}(zMN-\am)-\am\}\cr
&& 
+\sum_{n=1}^\infty \bm\otimes  e^{2z\ap}\frac {(-2z\ap)^n}{n!}
+z \sum_{n=1}^\infty \am\aa\otimes  e^{2z\ap}
\biggl( \frac {(-2z\ap)^n}{n!} -\frac {(-z\ap)^n}{n!}\biggr)\cr
&&  +z\sum_{n=0}^\infty\am\aa\otimes ( e^{2z\ap} -  e^{z\ap})
\frac {(-z\ap)^{n}}{n!}
-z\sum_{n=1}^\infty\am\otimes e^{z\ap}\frac {(-z\ap)^n}{n!}\aa \cr
&&=1\otimes \bm + \bm\otimes e^{2 z\ap}-z\am \otimes e^{z\ap}\aa
-z\aa\otimes \am\cr
&&+\bm\otimes e^{2 z\ap}(e^{-2 z\ap}-1) 
 +z\am\aa\otimes e^{2z\ap}(e^{-2z\ap}-e^{-z\ap}) \cr
&&+z\am\aa\otimes  ( e^{2z\ap} -  e^{z\ap})e^{-z\ap}
-z\am\otimes e^{z\ap}(e^{-z\ap}-1)\aa\cr
&&=1\otimes\bm + \bm\otimes 1 -z\aa\otimes \am - z \am \otimes
\aa.
\label{cp}
\eea
Now, we have to compute (\ref{cj}) with $f\equiv -z\ap\otimes\aa$ and
$g$ instead of $\Delta(X)$. In particular, 
\be
\begin{array}{l}
[-z\ap\otimes\aa,g]=z \{ e^{z\ap}(\am-zMN)+\am\}\otimes\aa
-z^2\bb e^{z\ap}\otimes\aa^2\cr
\qquad +2z\ap\otimes\bm +z^2\ap\otimes\am\aa-z(e^{z\ap}-1)\otimes
\am\aa-z^2\ap\aa\otimes \am ,\cr
\end{array} 
\label{cq}
\ee 
\be
\begin{array}{l}
[-z\ap\otimes\aa,[-z\ap\otimes\aa,g]]=
2z^2 \bb e^{z\ap}\otimes\aa^2+4z^2\ap^2\otimes\bm\cr
\qquad 
 +3z^3\ap^2\otimes\am\aa 
-2z^2\ap(e^{z\ap}-1)\otimes
\am\aa-z^3\ap^2\aa\otimes \am ,\cr
\end{array} 
\label{cr}
\ee 
and for $n\ge 3$:
\be
\begin{array}{l}
 [-z\ap\otimes\aa,[\dots [-z\ap\otimes\aa,g]]^{n)}\dots]=
z (2^n-1) (z\ap)^n\otimes\am\aa \cr
  + ( 2z \ap)^n\otimes\bm  - z n
(z\ap)^{n-1}(e^{z\ap}-1)\otimes
\am\aa -z(z\ap)^n \aa\otimes \am.\cr
\end{array} 
\label{cs}
\ee 
In this way, we complete the proof:
\bea
&&e^{-z\ap\otimes\aa}\,g \,e^{z\ap\otimes\aa}
=g+z \{ e^{z\ap}(\am-zMN)+\am\}\otimes\aa\cr
&&  +  
z\sum_{n=1}^\infty\biggl(  \frac {(2z\ap)^n}{n!}- 
\frac {(z\ap)^n}{n!}\biggr)\otimes\am\aa 
+\sum_{n=1}^\infty \frac {(2z\ap)^n}{n!}\otimes\bm\cr
&& - z \sum_{n=0}^\infty  \frac {(z\ap)^n}{n!}(e^{z\ap}-1)\otimes
\am\aa - z \sum_{n=1}^\infty\frac {(z\ap)^n}{n!}\aa\otimes\am\cr
&&=1\otimes\bm +\bm\otimes 1 -z\aa\otimes\am
+ze^{z\ap}(\am-zMN)\otimes\aa\cr
&& 
+z(e^{2z\ap}-e^{z\ap})\otimes\am\aa
+(e^{2z\ap}-1)\otimes\bm\cr
&& - z e^{z\ap} (e^{z\ap}-1)\otimes
\am\aa  - z (e^{z\ap}-1) \aa\otimes\am=\sigma\circ\Delta(\bm).
\eea

Finally, notice that ${\cal R}^{-1}=\sigma\circ{\cal
R}$, that is, $\cal R$ is a triangular $R$-matrix.


\end{document}